\documentclass[a4paper,12pt]{article}
\usepackage{color}
\usepackage{citesort}
\usepackage{amsmath,amssymb,latexsym,amsthm,graphicx,citesort}
\newtheorem{theorem}{Theorem}[section]
\newtheorem{lemma}[theorem]{Lemma}
\newtheorem{prop}[theorem]{Proposition}

\newcommand{\bproof}{{\bf Proof:~~}}
\newcommand{\eproof}{{\vrule height8pt width5pt depth0pt}\vspace{3mm}}

\newcommand{\mQ}{\mathcal{Q}}

\newcommand{\mC}{\mathcal{C}}

\newcommand{\mR}{\mathcal{R}}

\newcommand{\mD}{\mathcal{D}}

\newcommand{\mK}{\mathcal{K}}
\newcommand{\mL}{\mathcal{L}}
\newcommand{\mB}{\mathcal{B}}

\newcommand{\llg}{\lambda}

\newcommand{\bg}{\beta}

\newcommand{\og}{\omega}
\newcommand{\Og}{\Omega}
\newcommand{\pdh}{\partial}
\newcommand{\RR}{{\rm I\kern -1.6pt{\rm R}}}

\title{Range conditions for a spherical mean transform and global extendibility of solutions of Darboux equation}
\author{Mark Agranovsky and Linh V. Nguyen}

\begin{document}

\maketitle

\begin{abstract} We describe the range of the spherical Radon transform which evaluates integrals of a function
in $\RR^n$ over all spheres centered on a given sphere. Such transform attracts much attention
due to its applications in approximation theory and  (thermo- and photoacoustic) tomography.
Range descriptions for this transform have been obtained recently. They include two types of conditions:
orthogonality condition and, for even $n$, a moment condition. However, it was found out later that, in any dimensions, the moment condition follows from the orthogonality one, and therefore can be dropped. In terms of Darboux equation, which describes spherical means, it indirectly implies that solutions of certain boundary value problems in a domain automatically extend outside of the domain. In this article, we present a direct  proof of this global extendibility phenomenon for Darboux equation.
Correspondingly, it delivers an alternative proof of the range characterization theorem.

\end{abstract}

\section{Introduction}\label{S:Intro}
Let $f$ be a function defined on $\RR^n$. Consider the spherical mean Radon transform $$\mR f(x,t):= \frac{1}{\og_n} \int\limits_{S^{n-1}}f(x+t \theta) dA(\theta),$$ where $dA(\theta)$ is
the area measure on the unit sphere $S^{n-1} \subset \RR^n$  and $\og_n$ is the total measure of the unit sphere. This transform has been intensively investigated in the literature due to its applications in PDEs, medical imaging and geophysics (e.g., \cite{Jo,CH2,Be1,Be2,Ep,GGG,AKKun}).

It is well known that $G(x,t)=\mR f(x,t)$ satisfies the Darboux equation \cite{CH2,Jo}
\begin{eqnarray} \label{E:introd} \left\{\begin{array}{l} \mD G(x,t):=0, ~  (x,t) \in \RR^n \times \RR_+, \\ G(x,0) =f(x),~G_t(x,0) =0.\end{array} \right.\end{eqnarray} Here $\RR_+:=(0,\infty)$ and $\mD$ is the Darboux operator
\begin{equation}\label{N:Darboux} \mD:=\partial_t^2 + \frac{n-1}{t} \partial_t -\Delta_x.\end{equation}
Moreover, Aisgeirsson's theorem \cite[Ch.6]{CH2} states that any global $C^2$ solution $G(x,t)$ of (\ref{E:introd}) is the spherical means of its initial value: $G(x,t)=(\mR f)(x,t)$.

The problem of recovering a function $f$ of $n$ variables from the function $\mR f(x,t)$, which depends on $n+1$ parameters, is obviously overdetermined. Thus, it is natural to consider the restriction of $\mR f(x,t)$ to  an $n$-dimensional surface in the $(x,t)$-space $\RR^n \times \RR_+$. Different problems in analysis are related to different types of such surfaces.

Problems concerning the spherical mean operator restricted to cylindrical surfaces in $\RR^n \times \RR_+$ have attracted a lot of attention recently. By cylindrical surfaces we mean those of the form $\Gamma:=S \times \RR_+$ where $S$ is a hypersurface in $\RR^n$.
The problem of recovering finctions from their spherical means restricted to such surfaces has interesting applications in analysis, in particular, in approximation theory (see \cite{AQ,LP}). But perhaps more importantly, it serves as a mathematical model for thermoacoustic tomography, a novel medical imaging method (we refer the readers to \cite{FR,AKKun,KuKu} for detailed explanation). In that application, the restriction $g=\mR_S(f)$ of $\mR(f)$ to the set $\Gamma$ is the measured data, while $f$ is the image to be determined.

In this article, we are concerned with the problem of characterizing functions $g$ that belong to the range of the operator $\mR_S$.

Let us describe the problem in more precise terms (for the detailed exposition, we refer the reader to the articles
\cite{AmK06,FR06,AKQ}). First of all, in our considerations, the hypersurface  $S$ will be the unit sphere centered at the origin. We assume that $f \in C^{\infty}(\RR^n)$ is supported inside the closure $\overline{B}$ of the unit ball $B=\{|x| < 1\}$ (the class of such functions $f$ will be denoted by $C_0^{\infty}(\overline B)$).

We consider the problem  of characterizing  all functions $g(x,t)$ on the cylinder $\Gamma$ such that $g=\mR_S f$ for some function $f \in C_0^{\infty}(\overline B)$.

Some necessary conditions are almost obvious. First, the function $g(x,t)=\mR_S f(x,t)$ must be smooth on $\Gamma$. Second, $g(x,t)$ must vanish for all $t>2$ and also  vanish to infinite order at $t=0$. Thus, $g$ must satisfy
 {\it smoothness and support conditions}: $g \in C_0^{\infty}(S \times [0,2])$.

Another, less trivial, necessary condition can be derived from characterization of the spherical means by Darboux equation
(\ref{E:introd}). This condition, which we call {\it orthogonality condition} (see Theorem \ref{T:AKQ}),
is of Fredholm alternative type. It follows
from the existence of a solution of Darboux equation inside the cylinder $\Gamma$ (the fact that $S$ is a sphere is not important here).

The first complete range description for $n=2$ was obtained  in \cite{AmK06}.
It was proved that a function $g$  belongs to the image of $C_0^{\infty}(\overline B)$ under the operator $\mR_S$ if and only if it satisfies, besides the smoothness and support conditions and the above orthogonality condition,
an additional {\it moment condition}. Earlier, the necessity of the moment condition was observed in \cite{Patch}.

Further step in higher dimension was taken in \cite{FR06}. There a range characterization was obtained for odd  dimension and for a transform related to the wave equation, rather than Darboux equation.  The conditions in \cite{FR06} did not involve the moment conditions. A complete description of the range for the spherical mean transform in any dimension was obtained in  \cite{AKQ}.
As in \cite{AmK06}, the necessary and sufficient conditions in \cite{AKQ} fall into two groups: orthogonality  and moment conditions, although the moment conditions were needed for even dimensions only.

It had remained unclear whether the moment condition is really needed for even $n$, till the recent article \cite{AFK}. It was proved there that, regardless of parity of the dimension, the moment conditions follow from smoothness, support, and orthogonality conditions and therefore can be dropped.

The above result can be immediately translated to the language of Darboux equation. Namely, on one hand, it was proved in \cite{AKQ} that
the orthogonality condition for the data $g(x,t)$ is in fact the condition of existence of a solution $G=G^+$ of Darboux equation (\ref{E:introd}) inside the cylinder $\Gamma$ with the boundary data $G^+(x,t)=g(x,t)$ on $\Gamma$ and proper decay when $t \to + \infty$. On the other hand, spherical means are global solutions of Darboux equation. Therefore, $g$ is in the range of transform $\mR_S$ means that $G^+$ extends as a global solution to the entire space $\RR^n$.

The possibility of such an extension seems to be an interesting phenomenon by itself and
one may wish to have its direct proof.
In fact, such a proof was found in \cite{AKQ}, but in odd dimensions. In this article, we modify the construction of \cite{AKQ} to extend it to all dimensions. Correspondingly, we obtain an alternative  proof, universal for all dimensions,  of the range characterization theorem from \cite{AFK}.

\section{Main result}\label{S:main_result}
Let us recall here that $B$ is the unit ball centered at the origin, $S=\partial B$ is the unit sphere, and $\Gamma= S \times \RR_+$. We also denote by $J_{\mu}$ the Bessel function of the firts kind of order $\mu$ and $j_{\mu}(u)=u^{-\mu}J_{\mu}(u)$ the corresponding normalized Bessel function. The notation $C_0^{\infty}(\Omega)$ will stand for smooth functions in $\RR^n$  with the support in $\Omega$.

\subsection{ Formulation of main result}
The goal of this article is to present a direct  proof of the following result from \cite{AFK}:
\begin{theorem} \label{T:main}
Let $g$ be a function defined on the cylinder $\Gamma$. Then there exists $f \in C_0^{\infty} (\overline B)$ such that $g=\mR_S(f)$ if and only if the following conditions hold:
\begin{itemize}
\item[a)] {\bf Smoothness and support conditions:} $g \in C_0^\infty(\Gamma)$ and $g(x,t)=0$ when $t>2$.
\item[b)] {\bf Orthogonality condition:} Let $-\llg^2$ be an eigenvalue of the Dirichlet Laplacian on $B$ and $\varphi_\llg$  a corresponding eigenfunction. Then $$\int\limits_\Gamma g(x,t)\partial_{\nu_x} \varphi_{\lambda}(x) j_{\frac{n-2}{2}}(\lambda t) t^{n-1}dt dA(x)=0,$$ where $\nu_x$ is the outward normal to $S$ at $x$.
\end{itemize}
\end{theorem}

{\bf Remark:} Since $S$ is the unit sphere, the Dirichlet eigenfunctions $\varphi_{\lambda}$ can be written  in polar coordinates as follows:
$$\varphi_{\lambda}(r\theta)=j_{\frac{n+m-2}{2}}(\lambda r)Y_m(\theta),$$
where $Y_m$ is a spherical harmonics of degree $m$. The Dirichlet condition for $\varphi_\llg$ on the unit sphere $S$ requires that $j_{(n+m-2)/2}(\lambda)=0.$  Choosing $Y_m=Y{m,k}, k=1,\cdots,d(m),$ elements of the basis in the space of all harmonics of degree $m$, one can write the condition b) in the equivalent form:
\begin{itemize}
\item[$b^{\prime})$]
$\int\limits_\Gamma g(\theta,t)j_{\frac{n-2}{2}}(\lambda t)t^{n-1}Y_{m,k}(\theta)dtdA(\theta)=0.$
\end{itemize}
This can be rephrased as follows:
\begin{itemize}
\item[$b^{\prime\prime})$]
$\widehat{g}_{m,k}(\lambda)=0,$ for all zeros $\lambda$ of the Bessel function $j_{\frac{n+m-2}{2}}(\lambda)$, where
\end{itemize}
$$\widehat{g}_{m,k}(\lambda)=\int\limits_0 \limits^\infty g_{m,k}(t)j_{\frac{n-2}{2}}(\lambda \ t) t^{n-1}dt$$
is Fourier-Bessel transform of $g_{m,k}$ of the Fourier coefficient
$$g_{m,k}(t) = \int\limits_{S} g(y,t) Y_{m,k}(\theta) d A(\theta).$$

Theorem \ref{T:main} can be reformulated in terms of Darboux equation. Namely, since the spherical means $G=\mR f$ is the unique solution for the Darboux equation (\ref{E:introd}), Theorem \ref{T:main} is equivalent to:

\begin{theorem} \label{T:main2}
Let $g$ be a function defined on the cylinder $\Gamma$. Then the following statements are equivalent:
\begin{itemize}
\item[i)] There exists $f \in C_0^{\infty} (\overline B)$ such that the following problem has a solution:
\begin{eqnarray}\label{E:global} \left\{\begin{array}{l}\mD G (x,t)= 0, ~ (x,t) \in \RR^n \times \RR_+,\\ G(x,t) = g(x,t),
\ (x,t) \in \Gamma, \\ G(x,0) =f(x),~G_t(x,0) =0, \ x \in \RR^n. \end{array} \right.\end{eqnarray}
\item[ii)] The conditions a) and b) of Theorem \ref{T:main} hold.
\end{itemize}
\end{theorem}
The proof of the implication $i) \Rightarrow ii)$ is quite simple and can be found in \cite{AKQ}. In the rest of this article, we will prove the converse implication.
\subsection{Orthogonality condition and existence of the internal solution}
In this subsection, we recall a result from \cite{AKQ}, which is our starting point for the proof of Theorem \ref{T:main2}. This result interprets the orthogonality condition as the existence of some internal solution for the Darboux equation as follows:
\begin{theorem}\label{T:AKQ} Let $T>0.$ Consider the following backward boundary initial value problem for Darboux equation:
\begin{eqnarray*} \left\{\begin{array}{l} \mD G(x,t)=0, ~  (x,t) \in \Omega \times \RR_+, \\ G(x,t) =g(x,t), \ (x,t) \in \partial \Omega \times \RR_+, \\ G(x,t) =0, \ (x,t) \in B \times [T,\infty). \end{array} \right.\end{eqnarray*}
where $\Omega \subset \RR^n$ is a bounded domain with smooth boundary and
$g \in C^\infty_0(\partial \Omega \times [0,T])$. Then the following statements are equivalent:
\begin{itemize}
\item [a)] The  solution $G \in C^\infty(\overline{\Og} \times \RR_+)$ is smooth at $t=0$ and $G(x,t)=0$ for $t>T$.
\item [b)] The boundary data $g$ satisfies the orthogonality condition \begin{equation} \label{E:orth}\int\limits_{\partial \Omega}\int\limits_0\limits^{\infty}g(x,t)\partial_{\nu_x} \varphi_{\lambda}(x) j_{\frac{n-2}{2}}(\lambda t)t^{n-1}dA(x) dt=0,\end{equation} for all pairs of Dirichlet eigenvalue-eigenfunctions $\left(\llg, \varphi_{\lambda}\right)$.

\end{itemize}
\end{theorem}

The detailed proof can be found in \cite{AKQ}, and here we only briefly explain  its main idea. The implication
a) $\Rightarrow$ b) immediately follows from Stokes formula applied in $x$ variable and integration by parts with respect to $t.$
Vice versa, assume that the orthogonality condition (\ref{E:orth}) holds.
The equation is nonsingular for all $t>0$, and the unique solution exists in $B \times (0,T]$.
One needs only to show that $G$ is not singular at $t=0$ and $G_t(x,0)=0$. This is done by applying Stokes formula in the domain
$B \times (\epsilon,T], \epsilon \ge 0, $ using the orthogonality condition (\ref{E:orth}) and letting $\epsilon \to 0.$ .

{\bf Remark:} The above orthogonality condition can be understood as a Fredholm alternative statement for solvability of the  boundary value problem for the equation $\mD G=0$ in the solid cylinder $\Og \times \RR_+$. Indeed, denote by $Pg(x,t)$ the harmonic extension (Poisson integral) in $x$ of the boundary data $g(x,t), x \in S.$ By setting $H=G-Pg$, one can rewrite the boundary value problem \ref{E:introd} for $G$  as a boundary value problem for the nonhomogeneous equation
$$\mD H=-(\partial_t^2+\frac{n-1}{t}\partial_t)Pg,$$ with zero boundary data $H(x,t)=0, \ x \in S.$
Now, Fredholm alternative claims that for solvability the right hand side of the equation must be orthogonal to all solutions of the homogeneous adjoint equation. Since Darboux operator with Dirichlet  boundary conditions is self-adjoint, it suffices to check the orthogonality of the right hand side to separable solutions $u_{\lambda}(x,t)=j_{\lambda}(t)\varphi_{\lambda}(x)$, which constitute a complete system of  solutions. Then Fredholm alternative condition becomes
$$\int\limits_\Og \int\limits_{\RR_+}u_{\lambda}(x,t)Pg(x,t)t^{n-1}dtdx=0,$$
which is just our orthogonality condition (\ref{E:orth}), if uses the Stokes formula to replace the integration in $x$ over $\Og$  by the surface integration over $\partial \Og$.

\subsection{Theorem \ref{T:main2} as extendibility theorem for Darboux equation}
Let $\mC = B \times \RR_+, \Gamma=S \times \RR_+,$ and suppose that  the conditions $a)$ and $b)$ of Theorem \ref{T:main} are staisfied.
Consider the backward initial boundary value Darboux problem
\begin{eqnarray*} \left\{\begin{array}{l} \mD G(x,t)=0, ~  (x,t) \in \mC, \\ G(x,t) =g(x,t), \ (x,t) \in \Gamma, \\
G(x,t) =0,  \ (x,t) \in B \times [2,\infty). \end{array} \right.\end{eqnarray*}
Theorem \ref{T:AKQ} for the case $\Og =B$ and $T=2$ shows that the solution $G^+ \in C^\infty(\overline{\mC})$ of this problem is non-singular at $t=0$ and $G_t(x,0)=0.$ On the other hand, if a global solution $G$ for (\ref{E:global}) exists,
it also solves Darboux equation in the cylinder $\mC$ and satisfies the same boundary and initial conditions: $G(x,t)=g(x,t),
(x,t) \in \Gamma,$ and
$G(x,t)=0$ for $t>2$. Due to the uniqueness of the solution, $G$ and $G^+$ coincide on $\mC$. That is, the existence of $G$ is equivalent to the global extendibility for $G^+$.

The above argument shows that the implication $ii) \Rightarrow i)$ in Theorem \ref{T:main2} can be reduced to the following extendibility result:
\begin{theorem}\label{T:extend}
Let $G^+ \in C^{\infty}(\overline \mC)$ be a solution of Darboux equation $$\mD G^+(x,t)=0, \ (x,t) \in \mC,$$ such that
$G_t(x,0)=0$ and $G^+(x,t)=0$ for $ t>2$ and $x \in B$.
If the boundary value $g=G^+|_{\Gamma}$ belongs to $C_0^{\infty}(S \times [0,2]),$
then $G^+$  extends to a global solution $G$ of the Darboux equation (\ref{E:global}), with the $C^{\infty}$ initial value
$f^*(x)$ that is $G^+(x,0)$ extended by zero outside $B$.
\end{theorem}

\section {Proof of Theorem \ref{T:extend}}
The proof rests on the following two propositions. The first one says that if the initial data $f$ of the internal solution $G^+$ vanishes to infinite order on the unit sphere $S$, then the extendibility holds.

\begin{prop}\cite{AKQ} \label{P:p1}
Let $G^+ \in C^\infty(\overline{\mC})$ be a solution of Darboux equation $$\mD G^+(x,t)=0,~ \ (x,t) \in \mC$$
and $G^+(x,t)=0$ for all $ t>2$ and $ x \in B$. If $f(x):=G(x,0)$ vanishes to infinite order on the unit sphere $|x|=1,$ then $G^+$ extends to a global solution $G(x,t)$ of Darboux equation. This global solution is given by the spherical means $G(x,t) =(\mR f^*)(x,t)$ of the function $f^*$ in $\RR^n$
 obtained from $f$ by  the zero extension outside of the unit ball.
\end{prop}

The second proposition proves that the smoothness and support conditions for the boundary data $g$ imply required vanishing to infinite order of the function $f$ on the sphere $S.$

\begin{prop}\label{P:p2} Let $G^+ \in C^\infty(\overline{\mC})$ be a solution of Darboux equation $$\mD G(x,t)=0,~ \forall (x,t) \in \mC$$ such that $G(x,t)=0$ for  $t>2$ and $x \in B$. If the boundary value $g=G^+|_{\Gamma}$ belongs to $C_0^\infty(S \times [0,2])$  then $f(x)=G(x,0)$ vanishes to infinite order on the sphere $|x|=1$.
\end{prop}

We will prove Propositions \ref{P:p1} and \ref{P:p2} in separate sections. Meanwhile, we derive Theorem \ref{T:main} from these two propositions.

{\bf The proof of Theorem \ref{T:extend}}  is just simple combination of Propositions \ref{P:p1} and \ref{P:p2}.
Indeed, let $G^+(x,t)$ be the internal solution as in Theorem \ref{T:extend}. Then Proposition \ref{P:p2} implies that $f(x):=G^+(x,0)$ has zero of infinite order on the unit sphere $|x|=1$. By Proposition \ref{P:p1}, the (globally defined) spherical means $G(x,t):=\mR f^*(x,t)$ coincide with $G^+(x,t)$ inside the cylinder $\Gamma$ and therefore provide the extension of the internal solution $G^+.$ This proves Theorem \ref{T:extend} and hence  proves the equivalent Theorem \ref{T:main}.

\section{Proofs of the  propositions}

\subsection{Proposition \ref{P:p1}}
As we have already mentioned, Proposition \ref{P:p1} is  proven in \cite{AKQ}. We will present the proof here for the sake of completeness.

Since $f(x)=G^+(x,0)$ vanishes to infinite order on the boundary of $B$, the zero extended function $f^*$ belongs to $C^{\infty}(\RR^n)$.
Then the natural candidate for the extended solution is given by the spherical means of  $f^*$:
$$G = \mR f^*.$$
This function is globally defined and, since $f^*$ is smooth, belongs to $C^\infty(\RR^n \times \RR_+)$. It is  a global solution of Darboux equation and our goal is to prove that $G^+$ and $G$ coincide in the solid cylinder $\overline{\mC}:=\overline{B} \times [0,\infty)$.

First observe that both solutions share the same initial data at $t=2$:
$$G(x,2)=G^+(x,2)=G_t(x,2)=G^+_t(x,2)=0.$$
Then, due to domain of dependence theorem \cite[p.696]{CH2}, \begin{equation}\label{dde}
G(x,t)=G^+(x,t)=0 \mbox{ for all } (x,t) \in \mK^+,
\end{equation}
where $\mK^+$ is the upward characteristic cone
$$\mK^+:=\{(x,t) \in \overline{B} \times [0,2]: t-|x|\geq 1\}.$$
Moreover, $G(x,t)$ and $G^+(x,t)$ also share the initial data at $t=0$:
$$G(x,0)=G^+(x,0)=f(x), \\ G_t(x,0)=G^+_t(x,0)=0.$$
Therefore, again by the dependence domain theorem, they coincide in the downward characteristic cone
with the base $\overline{B} \times \{0\},$ \ $$\mK^-:=\{(x,t) \in \overline{B} \times [0,1]: |x|+t \leq 1\}.$$

Hence, the difference $U(x,t):=G(x,t)-G^+(x,t)$ vanishes in the union  $\mK =\mK^+ \cup \mK^-$ of the two cones.
Besides, since both  $G$ and $G^+$ vanish for $t>2$, so does $U$.

Since $U$ satisfies Darboux equation inside the cylinder $\mC$, its Fourier-Bessel transform $$\widehat U(x,\alpha)=\int\limits_0\limits^{\infty}U(x,t)j_{\frac{n-2}{2}}(\alpha t)t^{n-1}dt$$ satisfies Helmholtz equation:
\begin{eqnarray*} \Delta_x \widehat{U}(x,\alpha)=-\alpha^2\widehat U(x,\alpha),~\forall x \in B.\end{eqnarray*}
Hence, $\widehat{U}(x,\alpha)$ is real analytic with respect to  $x \in B$.
Also
$$U(x,t)=0, x \in \mK$$
and  $U(x,t)=0$ for $ t >2.$ The union $\mK \cup (B\times [2,\infty))$  contains the entire ray $\{(0,t): 0 \leq  t <\infty \}$
and hence after taking Fourier-Bessel transform one conludes that $\widehat U(0,\alpha)=0.$
Since $U$ is smooth, the same argument can be applied to $D^\bg_xU$ to obtain: $$D^{\beta}\widehat U(0,\alpha)= \widehat{D^{\beta}U}(0,\alpha)=0.$$
Thus, $\hat{U}(.,\llg)$ vanishes to infinite order at $x=0$. Since $\hat{U}(.,\llg)$ is real-analytic, one concludes
that $\widehat U(x,\alpha)=0, \ x \in B$  for all $\alpha.$
Taking inverse Fourier-Bessel transform, we obtain $U=0$ and therefore $G=G^+$ in $\mC$.
This  completes the proof of Proposition \ref{P:p1}.

\subsection{Proof of Proposition \ref{P:p2}}
We want to prove that $f:=G^+(.,0)$ vanishes on the sphere $S$ to infinite order:
$$D^{\beta}_x f(x)=0, \ |x|=1.$$
First of all, recall that $G(x,2)=G_t(x,2)=0$ implies $G^+=0$ in the upward characteristic cone $\mK^+$: $$G^+(x,t)=0 , \ |x| \leq t-1, 1 \leq t \leq 2.$$
In particular, $G^+$ vanishes at the vertex of the cone $\mK^+$:
$$G^+(0,1)=0.$$
Since $G^+$ is smooth in the solid cylinder $\overline{\mC}:= \overline B \times [0,\infty)$, the same conclusion holds for for all derivatives of $G^+$:
\begin{equation} \label{E:dG} D_t^{j}D_x^{\beta}G^+(0,1)=0.\end{equation}
Here $j=0,1,...$ and $\beta$ is an arbitrary multiindex.

Now we relate $G^+$ to the spherical means of the initial value $f(x)$. Albeit we cannot assert so far that $G^+=\mR f$ (which is in fact  our final goal), we can claim  that the two functions coincide at least in the downward characteristic cone $\mK^-$: $$\mK^-=\{|x|<1-t, \ 0 \leq t \leq 1\}.$$ Indeed, both $G^+$ and $\mR f$ solve the equation $\mD G=0$ on $\mC$ and share the same initial values on $\overline{B}$: $G(x,0)=f, G_t(x,0)=0$. The conclusion now follows from the domain of dependence argument.

Since $(0,1)$ is the vertex of $\mK^-$, due to (\ref{E:dG}), we conclude that
\begin{equation}\label{E:DR=0}
(D_t^{j}D_x^{\beta} \mR f)(0,1)=0.
\end{equation}
Now our aim is to derive from (\ref{E:DR=0}) that $f(x)$ vanishes on the sphere $|x|=1$ along with all derivatives.

To this end, first observe that condition (\ref{E:DR=0}) is invariant with respect to action of the orthogonal group $O(n)$ and
hence it holds for any term $f_{m,k}(r)r^mY_{m,k}(\theta)$ in Fourier series of $f.$
\begin{equation}\label{E:dec}
f(x)=f(r\theta)=\sum_{m=0}^{\infty}\sum_{k=1}^{d(m)}f_{m,k}(r)r^mY_{m,k}(\theta),
\end{equation} where $r=|x|, |\theta|=1$ and $Y_{m,k}, k=1,...,d(m)$ is the orthonormal basic in the space of all spherical harmonics
of degree $m.$
Due to smoothness of $f$ in the closed ball $\overline B$ , the series (\ref{E:dec}) converges
uniformly with all derivatives and hence it suffices to prove that each term vanishes on the unit sphere to infinite order.

Thus, we can assume that $f$ is just a single term: \begin{equation}\label{E:single}f(x)=f_{m}(r)P_m(x),\end{equation} where
$$P_m(x)=r^mY_{m,k}(\theta)$$
is a
spatial harmonic  of degree $m.$
To prove that $f$ vanishes to infinite order for $|x|=1$, it suffices to prove that all the derivatives $f_m^{(j)}$
vanihs at $t=1$:
$$f^{(j)}_{m}(1)=0, j=0,1,...$$ We will prove this by constructing a system of linear equations that these numbers satisfy.

\begin{lemma} \label{L:eq} The following identities hold
\begin{enumerate}
 \item For any $i \geq 0$, \begin{eqnarray} \label{E:L_m} (\frac {d^i}{dr^i} \mL_ m f_{m})(1)=0, \end{eqnarray} where $\mL_m$ is the following differential operator of order $m$:
\begin{eqnarray}\label{O:L_m} \mL_m  = \prod\limits_{s=1}\limits^m \left(\frac{1}{n+2(m-s)}r \frac{d}{dr} +1 \right).\end{eqnarray}
\item For any $l \geq 0$, \begin{eqnarray} \label{E:D_m} (\mQ_m)^l( f_{m})(1) = 0,\end{eqnarray} where $\mQ_m$ is the following differential operator of order $2$:
\begin{eqnarray}\label{O:D_m} \mQ_m = \partial^2_r+\frac{n+2m-1}{t} \partial_r.\end{eqnarray}
\end{enumerate} \end{lemma}

\bproof
It will be convenient to introduce the operator $$\pi_{m}(g)(r)=\sum_{k=1}^{d(m)} g_{m,k}(r)r^mY_k^m(\theta).$$
that projects onto spherical harmonics of degree $m$.
We observe that differentiation in $x$ and the transform $\mR$ commute: $$D_x^{\beta}\mR f(x,t)=\mR (D^{\beta}f) (x,t).$$
Since the spherical mean with the center at $x=0$ and  radius $t=1$ is exactly the projection onto order zero  harmonics (constants),
identity (\ref{E:DR=0}) now reads as $$\frac{d^i}{dr^i} \pi_{0}(D^{\beta}f)(1)=0$$
for all $i \in Z_+$ and all multiindices $\beta$.

The projection $\pi_0$ of the derivatives of $f$ of the form (\ref{E:single}) was computed in
 \cite[formula (2.10)]{Ep}:
$$\pi_0 (D^{\beta}f)=\left(D^{\beta}P_m\right) \left(\mL_{m} f_{m}\right),\ m=|\beta|=\beta_1+\cdots+\beta_n,$$ where the differential operator $\mL_m$ is defined in  (\ref{O:L_m}). Since $P_m$ is a polynomial of degree $m$, we can choose in this formula the multiindex $\beta$ so that $D^{\beta}P_m$ is a non-zero constant. By allowing the index $i$ to be arbitrary, we arrive at the identity (\ref{E:L_m}).

As for the identity (\ref{E:D_m}), it comes from the equation $\mD G^+=0$, which means
$$\mB G^+(x,t)=\Delta_xG^+(x,t),~ \forall (x,t) \in \mC,$$
where $\mB$ is the Bessel operator acting on $t-$variable $$\mB=\partial_t^2+\frac{n-1}{t}.$$ Iterating the above identity,
one obtains
\begin{equation} \label{E:interate} \mB^lG^+(x,t)=\Delta^lG^+(x,t),~\forall (x,t) \in \mC.\end{equation}
Since $G^+ \in C^\infty(\overline{C})$, the above equality holds up to the boundary $\Gamma$. In particular, since $G(x,t)=g(x,t)$ for all $(x,t) \in \Gamma$ and $G(x,0)=f(x)$, we have
\begin{equation} \label{E:interate} \mB^lg(x,0)=\Delta^lf(x),~\forall x \in S.\end{equation}
Since $g$ vanishes to infinite order at $t=0$, one concludes that for all $x \in S$
$$\Delta^l f(x)=0, l=0,1,...$$
Now, taking into account that the Laplacian $\Delta$ acts on $m^{th}$-harmonic term as the operator
$$\mQ_m:=\partial_r^2+\frac{n-1+2m}{r}\partial_r,$$ we arrive at \begin{equation}\label{add_eq} \mQ_m^lf_ {m}=0, l=0,1,...
\end{equation} \eproof

Let us apply Lemma \ref{L:eq} when indices $j$ and $l$ run independently from 0 to $m-1.$
We can write the differential operators in (\ref{E:L_m}) and (\ref{E:D_m}) in the form $$\frac{d^i}{dr^i}\mL_m= A_{i,0}(r)+A_{i,1}(r)\frac{d}{dr}+\cdots +A_{i,i+m}(r)\frac{d^{i+m}}{dr^{i+m}},$$ and $$\mD^l=B_{l,l}(r)\frac{d}{dr^l}+\cdots+B_{l,2l}(r) \frac{d^{2l}}{dr^{2l}}.$$
Consider the vector$$F:=(f_{m}(1),f_{m}^{\prime}(1), \cdots, f_{m}^{(2m-1)}(1)),$$ which consist
of the first $2m$ derivatives (including that of order $0$) of $f_{m}$ at the point $r=1$.

Let $i=0,1,\cdots,m-1$ and $l=0,1,\cdots,m-1$ in Lemma \ref{L:eq}, we conlcude that $F$ satisfies the following $2m \times 2m$ linear system:
\begin{eqnarray}\label{E:system}\left\{\begin{array}{l}
A_{0,0}F_0+ A_{0,1}F_1+ \cdots + A_{0,m}F_m= 0\\
A_{1,1}F_1+ A_{1,2}F_2+  \cdots \cdots \cdots+ A_{1,m+1}F_{m+1}=0\\
\cdots\cdots\cdots\\
A_{m-1,m-1}F_{m-1}+A_{m-1,m}F_m+ \cdots\cdots\cdots\cdots\cdots + A_{m-1,2m-1}F_{2m-1}=0\\
B_{0,0} F_0=0\\
B_{1,1}F_1+ B_{1,2}F_2=0\\
B_{2,2}F_2+ B_{2,3}F_3+ B_{2,4}F_4=0 \\
\cdots\cdots\cdots\\
B_{m-1,m-1}F_{m-1}+\cdots\cdots\cdots\cdots + B_{m-1,2m-2}F_{2m-2}=0.\end{array}
\right.
\end{eqnarray}
Here the matrix coefficients are $A_{i,j}=A_{i,j}(1), B_{i,j}=B_{i,j}(1)$.

\begin{lemma}\label{L:nondeg}
The linear $2m \times 2m-$ system (\ref{E:system}) is nondegenerate.
\end{lemma}
We will prove this lemma later.
Assuming that the lemma is proven, we can complete the proof of Proposition \ref{P:p2}. Since the system (\ref{E:system}) is nondegenerate, one concludes that the first $2m$ derivatives of $f_{m}$ vanish:
$$f_{m}^{(j)}(1)=F_{j}=0, 0 \leq j \leq 2m-1.$$
To obtain the vanishing of higher order derivatives, we will exploit higher values for the index $i$ in (\ref{E:L_m}). Choosing $i=m$ in (\ref{E:L_m}) results in shifting of the vector of the unknowns to the right: $$A_{m,0}F_m+\cdots +A_{m,2m-1}F_{2m-1}+A_{m,m}F_{2m}=0$$ which along with $F_0=\cdots=F_{2m-1}=0$ implies $F_{2m}=0$ because $A_{m,m} \neq 0.$
Then  the next choice $i=m+1$ leads to $F_{2m+1}=0.$ Proceeding this way by taking successively $j=m,m+1,m+2,\cdots$ one obtains $F_{\nu}=0$ for all $\nu\geq 0$.

This completes the proof of Proposition \ref{P:p2}. \eproof

{\bf Remark:} As it was mentioned earlier,  Proposition \ref{P:p2} was proved for odd $n$ in \cite{AKQ}. The proof used  Weyl and Poisson-Sonine integral transforms, applied to the solution $G^+(x,t)$ in $t$ and $x$ variables correspondingly.
To have control over the derivatives of $G^+(x,0)$ on the sphere $|x|=1$ one needs the inverse transforms to be local (differential) operators which is the case only in odd dimensions. That is why the proof in \cite{AKQ} did not generalize to even dimension.

{\bf Proof of Lemma \ref{L:nondeg}}

Let us denote by $A_i$ the $ith$ row of the matrix of the first $m$ equations in the system (\ref {E:system})
and by $B_i$ the $ith$ row from the second group of $m$ linear equations in (\ref{E:system}).

One observes that the vectors $A_i, i=0,\cdots,m-1$ are linearly independent, as $A_{i,j}=0$ for all $j>m+i$ and $A_{i,m+i} \neq 0$.

Now we will use induction with respect to the length of the system of the vectors. Namely,
we will show that on each step successive addition of vectors $B_0,\cdots,B_{m-1}$ to the set  $\{A^0,..,A^{m-1}\}$ does not violate the linear independence of the set obtained on the previous step. Then, in $m$ steps, we will obtain the linear independence of the entire system.

Thus, our inductive assumption is that the system
$$\mathcal{S}_p:=\{A_1,\cdots,A_{m-1},B_0,\cdots,B_{p-1}\} $$
is linearly independent for some $p \leq m-1$. Now we want to check that it remains linearly independent after adding the next vector $B_p$. In other words, the vector $B_p$ is linearly independent from the set $\mathcal S_p$.

To this end, it suffices to find a vector $v_p \in \RR^{2m}$ that is orthogonal to $\mathcal S_p$ but not to $B_p$:
\begin{eqnarray*}\left<A_i,v_p\right> &=& 0, i=0,\cdots,m-1,\\ \left<B_j,v_p\right> &=& 0, j=0,\cdots,p-1, \\ \left<B_p,v_p\right> &\neq& 0.\end{eqnarray*}

Indeed, we take the function $\Psi_p(r)=r^{-n-2p}$ and construct the vector of successive derivatives at $r=1$:
$$v_p=(\Psi_p (1),\cdots,\Psi_p^{(2m-1)}(1)).$$

Recall that $$\mL_m= \prod\limits_{s=1}\limits^m \left(\frac{1}{n+2(m-s)}r \frac{d}{dr} +1 \right).$$ Since all the first order differential operators in the above product commute, we can rewrite
$$\mL_m= \left[\prod\limits_{s \neq m-q}\limits^m \left(\frac{1}{n+2(m-s)}r \frac{d}{dr} +1 \right)\right] \left(\frac{1}{n+2p}r \frac{d}{dr} +1 \right).$$
A simple observation gives $$\left(\frac{1}{n+2p}r \frac{d}{dr} +1 \right) \Psi_p (r)=0,\forall r>0.$$
Therefore, $$\mL_m \Psi_p(r)=0,~\forall r>0.$$
This, in particular, implies
$$\left(\frac{d^i}{dr^i}\mL_m \Psi_p \right) (1)=0, \ i=0,1,\cdots,m-1.$$ By the definition of $A_i$, this is exactly equivalent to the first group of equations in (\ref{E:system}):
$$\left<A_i,v_p\right>=0, i=0,\cdots,m-1.$$
Here, as we defined, $v_p=(\Psi_p(1),\cdots, \Psi_p^{(1)}(1))$ is the vector of successive derivatives of $\Psi_p$ evaluated at $r=1$.

The second group of the equations in (\ref{E:system}) comes, by the construction, from the iteration of the differential operator
$$\mQ_m=\partial_r^2+\frac{n-1+2m}{r}\partial_r$$
evaluated at $r=1$.
Straightforward computation yields  \begin{eqnarray*} \left[\pdh_r^2 + \frac{n-1+2m}{r}\pdh_r \right]^{l}\Psi_p(r) =
\left\{\begin{array}{l}C_l r^{-n-2(p+l)},~ 0 \leq l \leq m-1-p, \\ 0,~ m-p \leq l \leq m-1.\end{array} \right. \end{eqnarray*} Here $C_l$ are nonzero constants.

Substituting $r=1$ and recollecting the original construction of the vectors $B_0,\cdots,B_{m-1}$ one is led to $$\left<B_0,v_p\right>=\cdots =\left<B_{p-1},v_p\right>=0, \ \left<B_p,v_p\right>=C_{m-1-p} \neq 0.$$ This completes the proof of Lemma \ref{L:nondeg} and thus finishes the proof of Theorem \ref{T:main}.
\eproof

\section*{Acknowledgments}

The first author is partially supported by the grant from ISF (Israel Science Foundation) 688/08. This work was done while he was on his sabbatical leave at Department of Mathematics, Texas A\&M University. He thanks ISF for support and the department for the hospitality and presenting excellent conditions for work.

The second author is partially supported by the NSF DMS grants 0604778 and 0715090 and the grant KUS-CI-016-04 from King Abdullah University of Science and Technology (KAUST). He thanks the NSF and KAUST for the support.

The authors thank Peter Kuchment for useful  discussions and valuable remarks.

\end{document}